\documentclass[]{article}
\pdfoutput=1
\usepackage{amsmath}
\usepackage{amssymb}
\usepackage{geometry}
\usepackage{color}
\begin{document}
\title{Distribution evaluation of hypergeometric series}
\author{Ming Hao Zhao}
\date{}
\maketitle
\abstract{We evaluate several classes of high weight hypergeometric series via Gamma, polylogarithm and elliptic integrals, mainly through distribution relations.}
\tableofcontents
\section{Preliminaries}

\noindent \textbf{Lemma 1.} Let $p,q,m,n\in \mathbb N, 0\le m<n, |z|\le1, A=\{a_1,\cdots,a_p\}, B=\{b_1,\cdots,b_q\}, D_{n,m}(c)=\{\frac{c+m}n,\cdots,\frac{c+m+n-1}n\}, D_{n,m}(C)=\cup_{c\in C}D_{n,m}(c)$, then\\
{\small $$_{np+1}F_{nq+n}(\{1\}\cup D_{n,m}(A);D_{n,m}(\{1\}\cup B);z)=m! n^{m(p-q-1)-1}z^{-\frac mn}\frac{\prod_{j=1}^q (b_j)_m}{\prod_{i=1}^p (a_i)_m}\sum_{k=0}^{n-1}e^{-\frac{2\pi ikm}n}\ _pF_q(A;B;n^{q-p+1}e^{\frac{2\pi ik}n}z^{\frac 1n})$$}\\
\noindent Proof. Check that RHS has coefficient 0 for all non-integer powers (elementary), and that both sides have same coefficient for all $z^n$ (Gauss multiplication of $\Gamma$).\\

We call this distributional operation \textcolor{blue}{DIST($n,m,z,A,B$)}. If $A, B$ are clear and RHS is summable using a set of formulas $X$, we abbreviate this operation as \textcolor{blue}{$X$, DIST($n,m,z$)} (mainly used in section 2). As an example, since $_3F_2$ form of $\sin^{-1}(z)^2$ is well known, one have $\text{DIST}\left(4,3,1,\{1,1,1\},\left\{\frac{3}{2},2\right\}\right)$ equivalent to:\\

{\footnotesize $$\, _6F_5\left(1,1,1,\frac{5}{4},\frac{3}{2},\frac{7}{4};\frac{9}{8},\frac{11}{8},\frac{13}{8},\frac{15}{8},2;1\right)=\frac{35}{64} \left(\pi ^2+4 \sin ^{-1}\left(\sqrt[4]{-1}\right)^2-4 \sinh ^{-1}(1)^2-4 \sinh ^{-1}\left(\sqrt[4]{-1}\right)^2\right)$$ $$=\frac{35}{64} \left(\pi ^2-4 \log ^2\left(\sqrt{2}+1\right)-\log ^2\left(\sqrt{2}-\sqrt{2 \sqrt{2}+2}+1\right)-\log ^2\left(\sqrt{2}+\sqrt{2 \sqrt{2}+2}+1\right)+8 \tan ^{-1}\left(\sqrt{\sqrt{2}-1}\right)^2\right)$$}\\

\noindent \textbf{Lemma 2.} We have:\\ $$_pF_q(a,a_2,\cdots,a_p;a-m,b_2,\cdots,b_q;z)$$ $$=\sum_{k=0}^m \binom{m}k \frac{z^k\prod_{j=2}^p (a_j)_k}{(a-m)_k \prod_{j=2}^q (b_j)_k}\ _{p-1}F_{q-1}(a_2+k,\cdots,a_p+k;b_2+k,\cdots,b_q+k;z)$$\\
$$_pF_q(1,a_2,\cdots,a_p;n+1,b_2,\cdots,b_q;z)$$ $$=\frac{(-1)^{(p-q)n}n! \prod_{j=2}^q (1-b_j)_n}{z^n \prod_{j=2}^p (1-a_j)_n}\left(\ _{p-1}F_{q-1}(a_2-n,\cdots,a_p-n;b_2-n,\cdots,b_q-n;z)-\sum_{k=0}^{n-1} \frac{z^k \prod_{j=2}^p (a_j-n)_k}{k! \prod_{j=2}^q (b_j-n)_k}\right)$$\\
$$\ _pF_q(a_1,\cdots,a_p;a_1+1,\cdots,a_n+1,b_{n+1}\cdots,b_q;z)$$ $$=\sum_{k=1}^n \prod_{j=1,j\not=k}^n \frac{a_j}{a_j-a_k}\ _{p-n+1}F_{q-n+1}(a_k,a_{n+1}\cdots,a_p;a_k+1,b_{n+1},\cdots,b_q;z)$$\\

\noindent Proof. First: In summand of the $_pF_q$ (with index $n$), decompose the polynomial $p(n)=\frac{\Gamma(n+a)}{\Gamma(n+a-m)}$ in terms of falling factorial i.e. $1, n, n(n-1),\cdots, n(n-1)\cdots (n-m+1)$, separate the sum into $m+1$ parts and reindex. Second: Eliminate initial terms and reindex. Third: Partial fractions (note that $a_1+1,\cdots,a_n+1$ can be generalized to arbitrary $a_1+k_1,\cdots,a_n+k_n$ with $k_i \in \mathbb Z$).\\

We name them as STIR (note the similarity to Stirling polynomials), INIT, PFD repsectively. One may derive a large class of new formulas from old ones depending on these operations. As an example, using Dougall $_5F_4$ ($L_3$ in Lemma 3 below) and INIT ($n=2$) we readily have:\\
$$\, _6F_5\left(\frac{1}{2},1,\frac{3}{2},\frac{3}{2},\frac{3}{2},\frac{9}{4};\frac{5}{4},2,2,2,3;1\right)=\frac{32}{5} \left(1-\frac{8}{\pi ^2}\right)$$\\

\section{Dougall Gamma DIST}

Now recall Gauss/Kummer (argument 1, $\frac12$)/Dixon/Watson/Whipple/Clausen formula \cite{Ask} for Gamma evaluation of $_2F_1$ and $_3F_2$ (we omit them). For higher weights we have Dougall formulas:\\

\noindent \textbf{Lemma 3.} We have:\\
$$L_1=\, _4F_3\left(a,\frac{a}{2}+1,b,c;\frac{a}{2},a-b+1,a-c+1;1\right)$$ $$=\frac{\Gamma \left(\frac{a+1}{2}\right) \Gamma (a-b+1) \Gamma (a-c+1) \Gamma \left(\frac{a+1}{2}-b-c\right)}{\Gamma (a+1) \Gamma \left(\frac{a+1}{2}-b\right) \Gamma \left(\frac{a+1}{2}-c\right) \Gamma (a-b-c+1)}$$\\
$$L_2=\, _4F_3\left(a,\frac{a}{2}+1,b,c;\frac{a}{2},a-b+1,a-c+1;-1\right)=\frac{\Gamma (a-b+1) \Gamma (a-c+1)}{\Gamma (a+1) \Gamma (a-b-c+1)}$$\\
$$L_3=\, _5F_4\left(\frac{a}{2}+1,a,b,c,d;\frac{a}{2},a-b+1,a-c+1,a-d+1;1\right)$$ $$=\frac{\Gamma (a-b+1) \Gamma (a-c+1) \Gamma (a-d+1) \Gamma (a-b-c-d+1)}{\Gamma (a+1) \Gamma (a-b-c+1) \Gamma (a-b-d+1) \Gamma (a-c-d+1)}$$\\
$$L_4=\, _6F_5\left(\frac{a}{2}+1,a,b,c,d,e;\frac{a}{2},a-b+1,a-c+1,a-d+1,a-e+1;-1\right)$$ $$=\frac{\Gamma (a-d+1) \Gamma (a-e+1)}{\Gamma (a+1) \Gamma (a-d-e+1)} \, _3F_2(a-b-c+1,d,e;a-b+1,a-c+1;1)$$\\

\noindent Proof. These are direct consequences of the finite Dougall $_7F_6$ summation formula; see \cite{Ask} for $L_2, L_3, L_4$. Now let $d\to \frac{a+1}2$ in $L_3$ to obtain $L_1$.\\

\noindent \textbf{Lemma 4.} We have:\\
$$K_1=\, _2F_1\left(a,\frac{a+1}{3};\frac{2 (a+1)}{3};e^{\pm \frac{\pi i}3}\right)=\frac{3^{\frac{a}{2}-1} e^{\frac{1}{6} \pi  a (\pm i)} \Gamma \left(\frac{a}{3}\right) \Gamma \left(\frac{2 (a+1)}{3}\right)}{\Gamma \left(\frac{2}{3}\right) \Gamma (a)}$$\\
$$K_2=\, _2F_1\left(a+\frac{1}{2},\frac{1}{2}-a;1-a;\frac{1}{2} \left(1-\sqrt{2}\right)\right)=\frac{\sqrt{\pi } 2^{2 a-\frac{1}{4}} \left(\sqrt{2}-1\right)^{a}  \Gamma (1-a)}{\Gamma \left(\frac{5}{8}-\frac{a}{2}\right) \Gamma \left(\frac{7}{8}-\frac{a}{2}\right)}$$\\
$$K_3=\, _3F_2\left(\frac{1}{2},\frac{1}{2}-a,a+\frac{1}{2};1-a,a+1;1\right)=\frac{\pi ^2 a \csc (\pi  a)}{\Gamma \left(\frac{3}{4}\right)^2 \Gamma \left(\frac{3}{4}-a\right) \Gamma \left(a+\frac{3}{4}\right)}$$\\
$$K_4=\, _3F_2\left(\frac{1}{2},\frac{1}{2}-a,a+\frac{1}{2};1-a,a+1;-1\right)=\frac{\pi ^2 a \csc (\pi  a)}{\sqrt{2} \Gamma \left(\frac{5}{8}-\frac{a}{2}\right) \Gamma \left(\frac{7}{8}-\frac{a}{2}\right) \Gamma \left(\frac{a}{2}+\frac{5}{8}\right) \Gamma \left(\frac{a}{2}+\frac{7}{8}\right)}$$\\

\noindent Proof. See \cite{Ebi} for proof of $K_1, K_2$ using nonlinear transformations of $_2F_1$ and limiting argument. Moreover one may establish the following identity:\\
{\footnotesize $$\, _3F_2\left(\frac{1}{2},\frac{1}{2}-a,a+\frac{1}{2};1-b,b+1;z\right)=\, _2F_1\left(\frac{1}{2}-a,a+\frac{1}{2};1-b;\frac{1-\sqrt{1-z}}{2} \right) \, _2F_1\left(\frac{1}{2}-a,a+\frac{1}{2};b+1;\frac{1-\sqrt{1-z}}{2}\right)$$}\\
By setting $z\to 4t(1-t)$ and checking that both sides satisfy the same order 3 ODE and initial conditions (Mathematica command `DifferentialRootReduce' should do it well). Now let $b\to a, z\to \pm1$ in above identity, use Kummer $\frac12$ and $K_2$ to obtain $K_3, K_4$ respectively.\\

\noindent \textbf{Proposition 1.} We have:\\

$$P11=\, _6F_5\left(\frac{a}{2}+1,a,1-b,b,1-c,c;\frac{a}{2},a-b+1,a+b,a-c+1,a+c;-1\right)$$ $$=\frac{\pi  2^{1-2 a} \Gamma (a-b+1) \Gamma (a+b) \Gamma (a-c+1) \Gamma (a+c)}{\Gamma (a) \Gamma (a+1) \Gamma \left(\frac{a}{2}-\frac{b}{2}-\frac{c}{2}+1\right) \Gamma \left(\frac{a}{2}+\frac{b}{2}-\frac{c}{2}+\frac{1}{2}\right) \Gamma \left(\frac{a}{2}+\frac{c}{2}-\frac{b}{2}+\frac{1}{2}\right) \Gamma \left(\frac{a}{2}+\frac{b}{2}+\frac{c}{2}\right)}$$\\
$$P12=\, _6F_5\left(\frac{a}{2}+1,a,2 b,a-2 c+1,2 a-2 b-2 c+1,c;\frac{a}{2},a-2 b+1,a-c+1,2 c,-a+2 b+2 c;-1\right)$$ $$=\frac{\sqrt{\pi } \Gamma \left(c+\frac{1}{2}\right) \Gamma (a-2 b+1) \Gamma (a-c+1) \Gamma (-a+2 b+2 c)}{\Gamma (a+1) \Gamma \left(b+\frac{1}{2}\right) \Gamma \left(-b+c+\frac{1}{2}\right) \Gamma (a-b-c+1) \Gamma (-a+b+2 c)}$$\\
$$P13=\, _6F_5\left(\frac{a}{2}+1,a,2 b,a-b-c+\frac{1}{2},2 c,-a+2 b+2 c;\frac{a}{2},a-2 b+1,a-2 c+1,2 a-2 b-2 c+1,b+c+\frac{1}{2};-1\right)$$ $$=\frac{\sqrt{\pi } \Gamma (a-2 b+1) \Gamma (a-2 c+1) \Gamma \left(b+c+\frac{1}{2}\right) \Gamma (a-b-c+1)}{\Gamma (a+1) \Gamma \left(b+\frac{1}{2}\right) \Gamma \left(c+\frac{1}{2}\right) \Gamma (a-b-2 c+1) \Gamma (a-2 b-c+1)}$$\\
$$P14=\, _6F_5\left(\frac{1}{2},a,a+\frac{3}{4},2 a-\frac{1}{2},2 a-b,b;a-\frac{1}{4},2 a,a+\frac{1}{2},2 a-b+\frac{1}{2},b+\frac{1}{2};-1\right)$$ $$=\frac{\sqrt{\pi } \Gamma \left(a+\frac{1}{2}\right)^2 \Gamma \left(b+\frac{1}{2}\right) \Gamma \left(2 a-b+\frac{1}{2}\right)}{\Gamma \left(2 a+\frac{1}{2}\right) \Gamma \left(\frac{b}{2}+\frac{1}{2}\right)^2 \Gamma \left(a-\frac{b}{2}+\frac{1}{2}\right)^2}$$\\
$$P15=\, _6F_5\left(\frac{a}{2}+1,a,a-2 b+1,a-b+\frac{1}{2},a-b+\frac{1}{2},b;\frac{a}{2},a-b+1,2 b,b+\frac{1}{2},b+\frac{1}{2};-1\right)$$ $$=\frac{\sqrt{\pi } \Gamma \left(b+\frac{1}{2}\right)^3 \Gamma (a-b+1)}{\Gamma (a+1) \Gamma \left(\frac{a}{2}-\frac{b}{2}+\frac{3}{4}\right)^2 \Gamma \left(-\frac{a}{2}+\frac{3 b}{2}+\frac{1}{4}\right)^2}$$\\

\noindent Proof. Choose appropriate parameters in $L_4$ then apply Dixon/Watson/Whipple/Clausen to sum RHS.\\

\noindent \textbf{Proposition 2.} We have:\\

$$P21=\, _4F_3\left(a,a+\frac{1}{2},b,b+\frac{1}{2};\frac{1}{2},a-b+\frac{1}{2},a-b+1;1\right)$$ $$=\frac{1}{2} \left(\frac{\Gamma (1-4 b)}{\Gamma (1-2 b) \Gamma (2 a-4 b+1)}+\frac{\sqrt{\pi } 4^{-a}}{\Gamma \left(a+\frac{1}{2}\right) \Gamma (a-2 b+1)}\right) \Gamma (2 a-2 b+1)$$\\
$$P22=\, _6F_5\left(\frac{a}{3}+\frac{1}{3},\frac{a}{3}+\frac{2}{3},a-\frac{1}{3},\frac{a}{3},a,a+\frac{1}{3};\frac{1}{3},\frac{2}{3},\frac{2 a}{3}+\frac{1}{3},\frac{2 a}{3}+\frac{2}{3},\frac{2 a}{3};-1\right)$$ $$=\sqrt{\pi } 8^{-a} \left(\frac{a \Gamma (2 a)}{\Gamma \left(\frac{a+1}{2}\right) \Gamma \left(\frac{3 a}{2}+1\right)}+\frac{2^{5 a+1} 3^{-\frac{3 a}{2}-1} \sin \left(\frac{1}{6} \pi  (4-3 a)\right) \Gamma \left(a+\frac{1}{2}\right)}{\Gamma \left(\frac{2}{3}\right) \Gamma \left(a+\frac{1}{3}\right)}\right)$$\\
$$P23=\, _6F_5\left(\frac{1}{4},\frac{3}{4},\frac{1}{4}-\frac{a}{2},\frac{3}{4}-\frac{a}{2},\frac{a}{2}+\frac{1}{4},\frac{a}{2}+\frac{3}{4};\frac{1}{2},\frac{1}{2}-\frac{a}{2},1-\frac{a}{2},\frac{a}{2}+\frac{1}{2},\frac{a}{2}+1;1\right)$$ $$=\frac{1}{4} \pi ^2 a \csc (\pi  a) \left(\frac{\sqrt{2}}{\Gamma \left(\frac{5}{8}-\frac{a}{2}\right) \Gamma \left(\frac{7}{8}-\frac{a}{2}\right) \Gamma \left(\frac{a}{2}+\frac{5}{8}\right) \Gamma \left(\frac{a}{2}+\frac{7}{8}\right)}+\frac{2}{\Gamma \left(\frac{3}{4}\right)^2 \Gamma \left(\frac{3}{4}-a\right) \Gamma \left(a+\frac{3}{4}\right)}\right)$$\\
$$P24=\, _6F_5\left(\frac{1}{4},\frac{3}{4},\frac{a}{2}+\frac{1}{2},a-\frac{1}{4},\frac{a}{2},a+\frac{1}{4};\frac{1}{2},\frac{a}{2}+\frac{1}{4},\frac{a}{2}+\frac{3}{4},a,a+\frac{1}{2};1\right)$$ $$=\frac{\pi  2^{-2 a-\frac{1}{2}} \Gamma \left(a+\frac{1}{2}\right)^2}{\Gamma \left(\frac{a}{2}+\frac{3}{8}\right)^2 \Gamma \left(\frac{a}{2}+\frac{5}{8}\right)^2}+\frac{\Gamma \left(\frac{1}{4}\right)^2 \Gamma \left(a+\frac{1}{2}\right)^2}{4 \pi  \Gamma \left(a+\frac{1}{4}\right)^2}$$\\
$$P25=\, _6F_5\left(\frac{3}{4},\frac{5}{4},\frac{a}{4}+\frac{3}{8},\frac{a}{4}+\frac{7}{8},\frac{a}{2}+\frac{1}{2},\frac{a}{2};\frac{1}{2},\frac{a}{4}+\frac{5}{8},\frac{a}{4}+\frac{9}{8},\frac{a}{2}+\frac{3}{4},\frac{a}{2}+\frac{5}{4};1\right)$$ $$=\frac{1}{2} \pi  \Gamma \left(\frac{a}{2}+\frac{5}{4}\right)^2 \left(\frac{2^{3 a-1} \left(\Gamma \left(\frac{a+1}{4}\right) \Gamma \left(\frac{a+2}{4}\right)-\Gamma \left(\frac{a}{4}\right) \Gamma \left(\frac{a+3}{4}\right)\right)^2}{\pi ^3 \Gamma (a)^2}+\frac{1}{\Gamma \left(\frac{5}{4}\right)^2 \Gamma \left(\frac{a+1}{2}\right)^2}\right)$$\\
$$P26=\, _7F_6\left(\frac{a}{4}+1,\frac{a}{2}+\frac{1}{2},\frac{a}{2},\frac{b}{2}+\frac{1}{2},\frac{b}{2},\frac{c}{2}+\frac{1}{2},\frac{c}{2};\frac{1}{2},\frac{a}{4},\frac{a}{2}-\frac{b}{2}+\frac{1}{2},\frac{a}{2}-\frac{b}{2}+1,\frac{a}{2}-\frac{c}{2}+\frac{1}{2},\frac{a}{2}-\frac{c}{2}+1;1\right)$$ $$=\frac{1}{2} \left(\frac{\Gamma (a-b+1) \Gamma (a-c+1)}{\Gamma (a+1) \Gamma (a-b-c+1)}+\frac{\Gamma \left(\frac{a+1}{2}\right) \Gamma (a-b+1) \Gamma (a-c+1) \Gamma \left(\frac{a+1}{2}-b-c\right)}{\Gamma (a+1) \Gamma \left(\frac{a+1}{2}-b\right) \Gamma \left(\frac{a+1}{2}-c\right) \Gamma (a-b-c+1)}\right)$$\\

\noindent Proof. P21: Gauss, Kummer $-1$, DIST(2,0,1). P22: Kummer $-1$, $K_1$ and its conjugate, DIST(3,0,$-1$). P23: $K_3, K_4$, DIST(2,0,1). P24, P25: Kummer $-1$, Clausen, DIST(2,0,1). P26: $L_1, L_2$, DIST(2,0,1).\\

\noindent \textbf{Proposition 3.} For $1\le k \le 6$ we have $_9F_8(Ak;Bk;1)=P3k$, where:\\

{\small $$A1=\left\{\frac{a}{4}+1,\frac{a}{2}+\frac{1}{2},\frac{a}{2},\frac{a}{2}-b+\frac{1}{2},\frac{a}{2}-b+1,\frac{3 a}{4}-b+\frac{1}{4},\frac{3 a}{4}-b+\frac{3}{4},\frac{b}{2}+\frac{1}{2},\frac{b}{2}\right\}$$ $$B1=\left\{\frac{1}{2},\frac{a}{4},\frac{a}{2}-\frac{b}{2}+\frac{1}{2},\frac{a}{2}-\frac{b}{2}+1,b,b+\frac{1}{2},-\frac{a}{4}+b+\frac{1}{4},-\frac{a}{4}+b+\frac{3}{4}\right\}$$
$$P31=\frac{4^{b-1} \Gamma \left(b+\frac{1}{2}\right) \Gamma (a-b+1) \left(\frac{\Gamma \left(\frac{a}{4}+\frac{1}{4}\right) \Gamma \left(-\frac{a}{4}+b+\frac{3}{4}\right)}{\Gamma \left(\frac{3 a}{4}-b+\frac{3}{4}\right) \Gamma \left(-\frac{3 a}{4}+2 b+\frac{1}{4}\right)}+\frac{\Gamma \left(-\frac{a}{2}+2 b+\frac{1}{2}\right) \Gamma \left(-\frac{3 a}{2}+3 b-\frac{1}{2}\right)}{\Gamma \left(-\frac{a}{2}+b+\frac{1}{2}\right) \Gamma \left(-\frac{3 a}{2}+4 b-\frac{1}{2}\right)}\right)}{\sqrt{\pi } \Gamma (a+1)}$$\\

$$A2=\left\{\frac{1}{4}-\frac{a}{4},\frac{3}{4}-\frac{a}{4},\frac{a}{4}+1,\frac{a}{2}+\frac{1}{2},\frac{a}{2},\frac{1}{2}-\frac{b}{2},1-\frac{b}{2},\frac{b}{2}+\frac{1}{2},\frac{b}{2}\right\}$$ $$B2=\left\{\frac{1}{2},\frac{3 a}{4}+\frac{1}{4},\frac{3 a}{4}+\frac{3}{4},\frac{a}{4},\frac{a}{2}-\frac{b}{2}+\frac{1}{2},\frac{a}{2}-\frac{b}{2}+1,\frac{a}{2}+\frac{b}{2},\frac{a}{2}+\frac{b}{2}+\frac{1}{2}\right\}$$
$$P32=\frac{\Gamma \left(\frac{3 a}{2}+\frac{1}{2}\right) \Gamma (a-b+1) \Gamma (a+b) \left(\frac{\Gamma \left(\frac{3 a}{2}-\frac{1}{2}\right)}{\Gamma (a) \Gamma \left(\frac{3 a}{2}-b+\frac{1}{2}\right) \Gamma \left(\frac{3 a}{2}+b-\frac{1}{2}\right)}+\frac{\pi ^{3/2} 2^{2-3 a}}{\Gamma \left(\frac{a}{2}\right) \Gamma \left(\frac{a}{4}-\frac{b}{2}+\frac{3}{4}\right) \Gamma \left(\frac{3 a}{4}-\frac{b}{2}+\frac{3}{4}\right) \Gamma \left(\frac{a}{4}+\frac{b}{2}+\frac{1}{4}\right) \Gamma \left(\frac{3 a}{4}+\frac{b}{2}+\frac{1}{4}\right)}\right)}{2 \Gamma (a+1)}$$\\

$$A3=\left\{\frac{a}{8}+\frac{1}{8},\frac{a}{8}+\frac{5}{8},\frac{a}{4}+1,\frac{a}{2}+\frac{1}{2},\frac{a}{2},\frac{3 a}{4}-\frac{b}{2}+\frac{1}{4},\frac{3 a}{4}-\frac{b}{2}+\frac{3}{4},\frac{b}{2}+\frac{1}{2},\frac{b}{2}\right\}$$ $$B3=\left\{\frac{1}{2},\frac{3 a}{8}+\frac{3}{8},\frac{3 a}{8}+\frac{7}{8},\frac{a}{4},\frac{a}{2}-\frac{b}{2}+\frac{1}{2},\frac{a}{2}-\frac{b}{2}+1,-\frac{a}{4}+\frac{b}{2}+\frac{1}{4},-\frac{a}{4}+\frac{b}{2}+\frac{3}{4}\right\}$$
$$P33=\frac{\Gamma (a-b+1) \Gamma \left(-\frac{a}{2}+b+\frac{1}{2}\right) \left(\frac{\sqrt{\pi } \Gamma \left(\frac{a}{4}+\frac{3}{4}\right) \Gamma \left(\frac{3 a}{4}+\frac{3}{4}\right)}{\Gamma \left(\frac{b}{2}+\frac{1}{2}\right) \Gamma \left(\frac{a}{4}-\frac{b}{2}+\frac{3}{4}\right) \Gamma \left(\frac{3 a}{4}-\frac{b}{2}+\frac{3}{4}\right) \Gamma \left(-\frac{a}{2}+\frac{b}{2}+\frac{1}{2}\right)}+\frac{\sec \left(\frac{1}{4} (3 \pi  a+\pi )\right) \cos \left(\frac{1}{4} \pi  (3 a-4 b+1)\right)}{\Gamma \left(\frac{1}{2}-\frac{a}{2}\right)}\right)}{2 \Gamma (a+1)}$$\\

$$A4=\left\{\frac{1}{4}-\frac{a}{4},\frac{3}{4}-\frac{a}{4},\frac{a}{4}+1,\frac{a}{2}+\frac{1}{2},\frac{a}{2},\frac{1}{2}-\frac{b}{2},1-\frac{b}{2},\frac{b}{2}+\frac{1}{2},\frac{b}{2}\right\}$$ $$B4=\left\{\frac{1}{2},\frac{3 a}{4}+\frac{1}{4},\frac{3 a}{4}+\frac{3}{4},\frac{a}{4},\frac{a}{2}-\frac{b}{2}+\frac{1}{2},\frac{a}{2}-\frac{b}{2}+1,\frac{a}{2}+\frac{b}{2},\frac{a}{2}+\frac{b}{2}+\frac{1}{2}\right\}$$
$$P34=\frac{\Gamma \left(\frac{3 a}{2}+\frac{1}{2}\right) \Gamma (a-b+1) \Gamma (a+b) \left(\frac{\pi  2^{1-2 a} \Gamma \left(\frac{a}{2}+\frac{1}{2}\right)}{\Gamma \left(\frac{a}{4}-\frac{b}{2}+\frac{3}{4}\right) \Gamma \left(\frac{3 a}{4}-\frac{b}{2}+\frac{3}{4}\right) \Gamma \left(\frac{a}{4}+\frac{b}{2}+\frac{1}{4}\right) \Gamma \left(\frac{3 a}{4}+\frac{b}{2}+\frac{1}{4}\right)}+\frac{\Gamma \left(\frac{3 a}{2}-\frac{1}{2}\right)}{\Gamma \left(\frac{3 a}{2}-b+\frac{1}{2}\right) \Gamma \left(\frac{3 a}{2}+b-\frac{1}{2}\right)}\right)}{2 \Gamma (a) \Gamma (a+1)}$$\\

$$A5=\left\{\frac{a}{4}+1,\frac{a}{2}+\frac{1}{2},\frac{a}{2},\frac{3 a}{8}-\frac{b}{4}+\frac{1}{8},\frac{3 a}{8}-\frac{b}{4}+\frac{5}{8},\frac{b}{2}+\frac{1}{2},-\frac{a}{4}+\frac{b}{2}+\frac{1}{4},-\frac{a}{4}+\frac{b}{2}+\frac{3}{4},\frac{b}{2}\right\}$$ $$B5=\left\{\frac{1}{2},\frac{a}{4},\frac{a}{2}-\frac{b}{2}+\frac{1}{2},\frac{a}{2}-\frac{b}{2}+1,\frac{3 a}{4}-\frac{b}{2}+\frac{1}{4},\frac{3 a}{4}-\frac{b}{2}+\frac{3}{4},\frac{a}{8}+\frac{b}{4}+\frac{3}{8},\frac{a}{8}+\frac{b}{4}+\frac{7}{8}\right\}$$
$$P35=\frac{2^{\frac{3 (a-1)}{2}-b} \Gamma \left(\frac{3 a}{4}-\frac{b}{2}+\frac{3}{4}\right) \Gamma \left(\frac{a}{4}+\frac{b}{2}+\frac{3}{4}\right)\left(\frac{\Gamma \left(\frac{a}{4}+\frac{1}{4}\right) \Gamma \left(\frac{a}{2}-\frac{b}{2}+1\right)}{\Gamma \left(\frac{b}{2}+\frac{1}{2}\right) \Gamma \left(\frac{3 a}{4}-b+\frac{3}{4}\right)}+\frac{\Gamma \left(\frac{3 a}{4}-\frac{3 b}{2}+\frac{1}{4}\right) \Gamma (a-b+1)}{\Gamma \left(\frac{3 a}{2}-2 b+\frac{1}{2}\right) \Gamma \left(\frac{a}{4}-\frac{b}{2}+\frac{3}{4}\right)}\right)}{\sqrt{\pi } \Gamma (a+1)}$$\\

$$A6=\left\{\frac{a}{8}+\frac{1}{8},\frac{a}{8}+\frac{5}{8},\frac{a}{4}+1,\frac{3 a}{8}+\frac{1}{8},\frac{3 a}{8}+\frac{1}{8},\frac{3 a}{8}+\frac{5}{8},\frac{3 a}{8}+\frac{5}{8},\frac{a}{2}+\frac{1}{2},\frac{a}{2}\right\}$$ $$B6=\left\{\frac{1}{2},\frac{a}{8}+\frac{3}{8},\frac{a}{8}+\frac{3}{8},\frac{a}{8}+\frac{7}{8},\frac{a}{8}+\frac{7}{8},\frac{3 a}{8}+\frac{3}{8},\frac{3 a}{8}+\frac{7}{8},\frac{a}{4}\right\},\ \ P36=\frac{\Gamma \left(\frac{a+3}{4}\right)^2 \left(\frac{\sqrt{\pi } \Gamma \left(\frac{a}{4}+\frac{3}{4}\right) \Gamma \left(\frac{3 a}{4}+\frac{3}{4}\right)}{\Gamma \left(\frac{5}{8}-\frac{a}{8}\right)^2 \Gamma \left(\frac{3 a}{8}+\frac{5}{8}\right)^2}+\frac{\sec \left(\frac{1}{4} (3 \pi  a+\pi )\right)}{\Gamma \left(\frac{1}{2}-\frac{a}{2}\right)}\right)}{2 \Gamma (a+1)}$$}\\

\noindent Proof. $L_3$,  certain $_5F_4$ formulas induced by $L_4$ (let $e\to \frac{a+1}2$ and choose appropriate $b,c,d$ to make RHS summable using Dixon/Watson/Whipple/Clausen), DIST(2,0,1).\\

\noindent Note 1: One may try DIST(2,1,1), DIST(3,1,$-1$), etc in various steps to derive analogous formulas.\\

\noindent Note 2: One may establish the FL expansion $(x(1-x))^{s-1}=B(s,s)\sum_{n=0}^\infty  \frac{(5/4)_n(1-s)_n(1/2)_n}{(1/4)_n(1/2+s)_n(1)_n} P_{2n}(2x-1)$ \cite{ZMH} using Dixon formula. Thus FL Parseval on $(x(1-x))^{s-1}, (x(1-x))^{t-1}$ gives\\
$$\, _5F_4\left(\frac{1}{2},\frac{1}{2},\frac{5}{4},1-s,1-t;\frac{1}{4},s+\frac{1}{2},t+\frac{1}{2},1;1\right)=\frac{B(s+t-1,s+t-1)}{B(s,s) B(t,t)}$$\\
Which is a symmetric special case of $L_3$.\\

\noindent Note 3: It is not easy to derive irreducible Gamma closed-forms of higher order series (e.g. $_{10}F_9, \cdots$). For instance calculating $P22(a\to \frac14), P24(a\to \frac16)$, DIST(2,0,1) yields a seemingly nontrivial result:\\
$$\, _{12}F_{11}\left(-\frac{1}{24},\frac{1}{24},\frac{1}{8},\frac{5}{24},\frac{7}{24},\frac{3}{8},\frac{11}{24},\frac{13}{24},\frac{5}{8},\frac{17}{24},\frac{19}{24},\frac{7}{8};\frac{1}{12},\frac{1}{6},\frac{1}{4},\frac{1}{3},\frac{5}{12},\frac{1}{2},\frac{7}{12},\frac{2}{3},\frac{3}{4},\frac{5}{6},\frac{11}{12};1\right)$$ $$=\frac{2 \sqrt{\sqrt{2}+1}+\frac{\sqrt[8]{3} \left(\sqrt{3}+1\right)  \Gamma \left(\frac{1}{3}\right) \Gamma \left(\frac{5}{12}\right)}{\Gamma \left(\frac{1}{4}\right)}\sqrt{\frac{\sqrt{\sqrt{3}+2}+2}{\pi }}}{12 \sqrt{2}}+\frac{\pi  \left(\frac{\Gamma \left(\frac{1}{4}\right)^2}{ \Gamma \left(\frac{5}{12}\right)^2}+\frac{\left(-\sqrt{2}+\sqrt{6}+4\right) \Gamma \left(\frac{7}{24}\right)^2}{ 2^{11/6} \Gamma \left(\frac{11}{24}\right)^2}\right)}{6 \Gamma \left(\frac{1}{3}\right)^2}$$\\
Which is in fact trivial, since DIST($12,0,1,\{-\frac {1}2\},\{\}$) evaluates the same series to an algebraic number:\\
$$\frac{1}{12} \sqrt{1-\sqrt[6]{-1}}-\frac{1}{12} (-1)^{5/6}+\frac{\sqrt[6]{-1}}{12}-\frac{1}{12} (-1)^{11/12} \sqrt{\sqrt[6]{-1}-1}+\frac{1}{12} \sqrt{1-(-1)^{5/6}}-\frac{1}{12} (-1)^{7/12} \sqrt{(-1)^{5/6}-1}$$ $$+\frac{\sqrt[8]{-1}}{6\ 2^{3/4}}-\frac{(-1)^{7/8}}{6\ 2^{3/4}}+\frac{\sqrt[12]{-1}}{4\ 3^{3/4}}-\frac{(-1)^{11/12}}{4\ 3^{3/4}}+\frac{1}{6 \sqrt{2}}$$\\
Transforming the first form into the second (i.e. eliminating Gammas) directly, as well as reducing the second to real radicals, are not so trivial though. See `Examples' section for more.\\

\section{Polylog, elliptic K/E DIST}

\subsection{Denominator 2}
Evidently all $_1F_0$ are trivial (see Note 3 in section above for an example). Moreover, all $_2F_1$ with half-integer parameters are expressible via elementary functions or elliptic-$K/E$.\\

Now denote $O=\mathbb Z+\frac12, E=\mathbb Z$. Thus one may separate all $_3F_2$ with half-integer parameters into 12 classes: $EEE|EE, EEE|EO, EEE|OO, EEO|EE,\cdots, OOO|OO$. For instance $EEE|OO$ denotes all series of form $_3F_2(i,j,k;m+\frac12,n+\frac12;z)$ with $i,j,\cdots$ integers. Using contiguous relations it is clear that 10 classes other than $OOO|OE, EEE|OO$ have polylog or elliptic-$K/E$ closed-forms. Among 10 classes, functions belong to $OOO|OO, EEO|OO, EOO|EE, OOO|EE$ are of considerable complexity (polylog for former two and elliptic for latter two). Here are examples of these 4 cases:\\
$$\sqrt{z} \, _3F_2\left(\frac{1}{2},\frac{1}{2},\frac{1}{2};\frac{3}{2},\frac{3}{2};-z\right)=\text{Li}_2\left(-\sqrt{z}-\sqrt{z+1}\right)-\text{Li}_2\left(-\sqrt{z}-\sqrt{z+1}+1\right)$$ $$-\frac{1}{2} \log ^2\left(\sqrt{z}+\sqrt{z+1}\right)+\log \left(\sqrt{z}+\sqrt{z+1}+1\right) \log \left(\sqrt{z}+\sqrt{z+1}\right)+\frac{\pi ^2}{12}$$\\
$$\sqrt{z} \, _3F_2\left(1,1,\frac{1}{2};\frac{3}{2},\frac{3}{2};-z\right)=-\text{Li}_2\left(\frac{1}{\sqrt{z}+\sqrt{z+1}}\right)+\text{Li}_2\left(-\frac{1}{\sqrt{z}+\sqrt{z+1}}\right)$$ $$+\log \left(\sqrt{z}+\sqrt{z+1}\right) \log \left(\frac{\sqrt{z}+\sqrt{z+1}-1}{\sqrt{z}+\sqrt{z+1}+1}\right)+\frac{\pi ^2}{4}$$\\
$$9 \pi  z^2 \, _3F_2\left(\frac{1}{2},\frac{3}{2},1;3,3;z\right)=-512 z K(z)+512 K(z)+144 \pi  z-128 z E(z)-896 E(z)+192 \pi$$\\
$$\pi ^2 \sqrt{1-z} \, _3F_2\left(\frac{1}{2},\frac{1}{2},\frac{3}{2};1,1;z\right)=8 K\left(\frac{1}{2}-\frac{\sqrt{1-z}}{2}\right) E\left(\frac{1}{2}-\frac{\sqrt{1-z}}{2}\right)-4 K\left(\frac{1}{2}-\frac{\sqrt{1-z}}{2}\right)^2$$\\
Note while calculating DIST relations, the last class $OOO|EE$ gives $K,E$ of argument $\frac12$ and $\frac{1}{2}-\frac{\sqrt{2}}{2}$, both of which are related to ESVs (elliptic singular values) $K(k_1), K(k_2)$ and quadratic transforms, say\\
$$K\left(\frac{1}{2}-\frac{\sqrt{2}}{2}\right)=\frac{\Gamma \left(\frac{1}{8}\right) \Gamma \left(\frac{3}{8}\right)}{2^{11/4} \sqrt{\pi }},E\left(\frac{1}{2}-\frac{\sqrt{2}}{2}\right)=\frac{\Gamma \left(\frac{5}{8}\right) \Gamma \left(\frac{7}{8}\right)+\frac{1}{8} \left(\sqrt{2}+1\right) \Gamma \left(\frac{1}{8}\right) \Gamma \left(\frac{3}{8}\right)}{2^{5/4} \sqrt{\pi }}$$\\
Thus all $_pF_q$ obtained by taking DIST on $OOO|EE$ are expressible via Gamma functions, in a way completely different from Dougall summation. As one may see from examples below, DIST on other classes also give interesting patterns of polylogarithm.\\

For $_4F_3$, $EEEE|EEO, EEEO|EEE, OOOO|OOO$ are the only classes not decomposable to $_3F_2$ that have polylog closed forms (they are obtained by brute-force integration, other than elementary means such as STIR/INIT/PFD). Functions in $EOOO|EEE$ are given closed-forms by INIT while that of several other classes are given by PFD. Examples of 3 irreducible cases:\\
$$\, _4F_3\left(\frac{1}{2},\frac{1}{2},\frac{1}{2},\frac{1}{2};\frac{3}{2},\frac{3}{2},\frac{3}{2};z\right)=\frac{\sin ^{-1}\left(\sqrt{z}\right)^3}{6 \sqrt{z}}+\frac{i \left(\text{Li}_3\left(2 z-2 i \sqrt{z-z^2}\right)-\text{Li}_3\left(2 z+2 i \sqrt{z-z^2}\right)\right)}{4 \sqrt{z}}$$\\
$$z \, _4F_3\left(\frac{1}{2},1,1,1;2,2,2;z\right)=-4 \text{Li}_2\left(\frac{1}{2}-\frac{\sqrt{1-z}}{2}\right)-8 \sqrt{1-z}+2 \log ^2\left(\sqrt{1-z}+1\right)$$ $$-4 \log (2) \log \left(\sqrt{1-z}+1\right)+8 \log \left(\sqrt{1-z}+1\right)+8+2 \log ^2(2)-8 \log (2)$$\\
$$z \, _4F_3\left(1,1,1,1;\frac{3}{2},2,2;z\right)=\text{Li}_3\left(-2 z-2 i \sqrt{1-z} \sqrt{z}+1\right)+2 i \text{Li}_2\left(-2 z-2 i \sqrt{1-z} \sqrt{z}+1\right) \sin ^{-1}\left(\sqrt{z}\right)$$ $$+\frac{2}{3} i \sin ^{-1}\left(\sqrt{z}\right)^3+2 \log (2) \sin ^{-1}\left(\sqrt{z}\right)^2+2 \log \left(\left(\sqrt{z}+i \sqrt{1-z}\right) \sqrt{z}\right) \sin ^{-1}\left(\sqrt{z}\right)^2-\zeta (3)$$\\

For $_5F_4$, $EEEEO|EEEE$ are expressible using complicated polyloarithms via integrating $EEEO|EEE$. See \cite{ZMH1}, Prop. 14 for an example featuring $\, _5F_4\left(1,1,1,1,\frac{3}{2};2,2,2,2;z\right)$ (we omit it for sake of simplicity). One may also obtain closed-forms of $EEEEE|EEEO$ via integration but they are even more complicated (the only one we shall use is $\, _5F_4\left(1,1,1,1,1;\frac{3}{2},2,2,2;z\right)$, which is due to \cite{5f4}). We are not aware of non-elementary closed-forms of higher order ones.\\

\subsection{Denominator 4}
Denote $P=\mathbb Z+\frac14, Q=\mathbb Z+\frac34$. Then all $_2F_1$ of class $PQ|O$ are expressible via elementary functions. By integrating functions of this type one may obtain polylog closed forms of class $E_k PQ|E_k O, O_kPQ|O_{k+1}$ ($E_k$ denotes $E\cdots E$ for $k$ times) or even $PPQ|PO, PQQ|QO$. All $_3F_2$ of class $EEO|PQ$ also admit elementary closed-forms. Several more classes like $OPQ|EE$ are expressible via elliptic K/E but they are of less interest to us. For denominator 4 case we shall not elaborate examples here but refer to the section below.\\

Another important result is that all $\, _2F_1(1,1;k;z), \, _2F_1(1,k-1;k;z), \, _2F_1(k-1,k-1;k;z)$ with $k\in \mathbb Q$ are expressible via elementary functions, due to Pfaff/Euler transforms (resulting in middle form), $_2F_1$ form of $\log(1-z)$ and DIST relation on $\sum_{n=0}^\infty \frac{z^n}{pn+q}$. When $k\in P, Q$ the expression is simple enough to carry out further integration, yielding closed-forms of some $PPP|PP, QQQ|QQ$, etc.\\

One may plug in special values or apply DIST relations in functions above to obtain hypergeometric evaluations. One may also multiply these functions by Beta and log kernel $x^a (1-x)^b, x^a \log^n(x)$ and integrate on $(0,1)$ to get more. For low order Pfaff/Euler/Thomae transforms may be applied again. Also STIR is always usable while generating high order $_pF_q$.

\clearpage
\section{Examples}
\subsection{Dougall Gamma DIST}
Choosing suitable $a,b,c,\cdots$ in DIST formulas of Prop. 1,2,3 gives the following.\\

\noindent Prop. 1:\\
$$\, _6F_5\left(\frac{1}{4},\frac{1}{4},\frac{1}{4},\frac{3}{4},\frac{3}{4},\frac{9}{8};\frac{1}{8},\frac{1}{2},\frac{1}{2},1,1;-1\right)=\frac{\sqrt[4]{2} \left(\sqrt{2}+1\right) \Gamma \left(\frac{3}{8}\right)^2}{\pi ^{3/2} \Gamma \left(\frac{1}{4}\right)}$$\\
$$\, _6F_5\left(-\frac{1}{4},\frac{1}{12},\frac{1}{8},\frac{1}{6},\frac{1}{2},\frac{7}{8};-\frac{1}{8},\frac{1}{4},\frac{7}{12},\frac{5}{8},\frac{2}{3};-1\right)=\frac{\left(\sqrt{6}+3\right) \Gamma \left(\frac{1}{4}\right) \Gamma \left(\frac{5}{12}\right) \Gamma \left(\frac{11}{24}\right)^2}{6 \sqrt{\pi } \Gamma \left(\frac{1}{3}\right) \Gamma \left(\frac{3}{8}\right)^2}$$\\
$$\, _6F_5\left(\frac{1}{8},\frac{1}{6},\frac{11}{24},\frac{11}{24},\frac{19}{24},\frac{17}{16};\frac{1}{16},\frac{1}{3},\frac{2}{3},\frac{2}{3},\frac{23}{24};-1\right)=\frac{32 \left(\sqrt{2-\sqrt{\sqrt{3}+2}}+2\right) \pi ^{5/2} \Gamma \left(\frac{13}{48}\right)^2}{3 \sqrt{6-3 \sqrt{\sqrt{3}+2}} \Gamma \left(\frac{1}{24}\right) \Gamma \left(\frac{1}{8}\right) \Gamma \left(\frac{1}{3}\right)^3 \Gamma \left(\frac{7}{16}\right)^2}$$\\
Prop. 2:\\
{\small $$\, _7F_6\left(\frac{1}{6},\frac{1}{6},\frac{1}{4},\frac{2}{3},\frac{2}{3},\frac{3}{4},\frac{9}{8};\frac{1}{8},\frac{1}{2},\frac{7}{12},\frac{7}{12},\frac{13}{12},\frac{13}{12};1\right)=\frac{\Gamma \left(\frac{1}{3}\right)^6 \left(\sqrt{3}+\frac{3 \sqrt{\sqrt{3}+2} \Gamma \left(\frac{1}{3}\right)^2}{\Gamma \left(\frac{1}{4}\right) \Gamma \left(\frac{5}{12}\right)}\right)}{48 \pi ^3}$$\\
$$\, _6F_5\left(\frac{1}{6},\frac{1}{4},\frac{1}{3},\frac{2}{3},\frac{3}{4},\frac{5}{6};\frac{5}{12},\frac{1}{2},\frac{7}{12},\frac{11}{12},\frac{13}{12};1\right)=\frac{4 \Gamma \left(\frac{1}{12}\right) \Gamma \left(\frac{5}{12}\right) \Gamma \left(\frac{1}{4}\right)^2+\sqrt{2} \Gamma \left(\frac{1}{24}\right) \Gamma \left(\frac{5}{24}\right) \Gamma \left(\frac{7}{24}\right) \Gamma \left(\frac{11}{24}\right)}{192 \pi ^2}$$\\
$$\, _6F_5\left(\frac{1}{6},\frac{11}{24},\frac{2}{3},\frac{3}{4},\frac{23}{24},\frac{5}{4};\frac{1}{2},\frac{17}{24},\frac{11}{12},\frac{29}{24},\frac{17}{12};1\right)=\frac{25\left(-3 \sqrt{3}+6+\frac{2^{2/3} \pi  \Gamma \left(\frac{1}{12}\right)^2 \Gamma \left(\frac{5}{12}\right)^2}{\Gamma \left(\frac{1}{3}\right)^6}-\frac{3\ 2^{5/6} \Gamma \left(\frac{5}{12}\right)^2}{\sqrt{\pi } \Gamma \left(\frac{1}{3}\right)}+\frac{9 \Gamma \left(\frac{1}{3}\right)^2 \Gamma \left(\frac{5}{12}\right)^2}{\pi  \Gamma \left(\frac{1}{4}\right)^2}\right)}{216} $$}\\
Prop. 3:\\
$$\, _9F_8\left(-\frac{1}{16},\frac{1}{8},\frac{1}{8},\frac{1}{4},\frac{7}{16},\frac{5}{8},\frac{5}{8},\frac{3}{4},\frac{17}{16};\frac{1}{16},\frac{3}{8},\frac{1}{2},\frac{1}{2},\frac{11}{16},\frac{7}{8},1,\frac{19}{16};1\right)=\frac{6 \sqrt{2} \left(\sqrt{\left(3-2 \sqrt{2}\right) \pi }+\frac{\sqrt{\pi } \Gamma \left(\frac{3}{16}\right) \Gamma \left(\frac{5}{16}\right)}{\Gamma \left(\frac{1}{16}\right) \Gamma \left(\frac{7}{16}\right)}\right)}{\Gamma \left(\frac{1}{4}\right)^2}$$\\
$$\, _9F_8\left(-\frac{1}{4},-\frac{1}{16},-\frac{1}{16},\frac{1}{16},\frac{1}{4},\frac{7}{16},\frac{7}{16},\frac{9}{16},\frac{7}{8};-\frac{1}{8},\frac{3}{16},\frac{5}{16},\frac{5}{16},\frac{1}{2},\frac{11}{16},\frac{13}{16},\frac{13}{16};1\right)$$$$=\frac{\pi  \left(\sqrt{2-\sqrt{2}}+2\right) \Gamma \left(\frac{5}{16}\right)^2+2 \sqrt{2 \pi } \Gamma \left(\frac{1}{4}\right) \Gamma \left(\frac{7}{16}\right)^2}{\left(\sqrt{2}+2\right)^{3/2} \Gamma \left(\frac{3}{8}\right)^2 \Gamma \left(\frac{7}{16}\right)^2}$$\\
$$\, _9F_8\left(\frac{1}{24},\frac{1}{12},\frac{5}{12},\frac{5}{12},\frac{13}{24},\frac{7}{12},\frac{11}{12},\frac{11}{12},\frac{29}{24};\frac{5}{24},\frac{1}{2},\frac{1}{2},\frac{5}{6},\frac{7}{8},1,\frac{4}{3},\frac{11}{8};1\right)$$$$=\frac{3 \sqrt{3} \Gamma \left(\frac{1}{3}\right) \left(\frac{3 \sqrt[6]{2} \Gamma \left(\frac{1}{12}\right) \Gamma \left(\frac{1}{3}\right)^2 \Gamma \left(\frac{5}{12}\right)^2}{\pi  \Gamma \left(\frac{1}{4}\right)}-\frac{\sqrt{2-\sqrt{2}} \left(\sqrt{3}-3\right) \Gamma \left(\frac{1}{24}\right) \Gamma \left(\frac{1}{8}\right) \Gamma \left(\frac{11}{24}\right)}{\Gamma \left(\frac{7}{24}\right)}\right)}{70\ 2^{5/6} \pi  \Gamma \left(\frac{1}{4}\right) \Gamma \left(\frac{5}{12}\right)}$$\\
Bonus: Similar to Note 3, choosing appropriate parameters in analogy of $P22$ (using DIST(3,1,$-1$) instead) and original $P24$, DIST(2,0,1) gives\\
{\small  $$S=\, _{12}F_{11}\left(\frac{1}{8},\frac{5}{24},\frac{7}{24},\frac{3}{8},\frac{11}{24},\frac{13}{24},\frac{5}{8},\frac{17}{24},\frac{19}{24},\frac{7}{8},\frac{23}{24},\frac{25}{24};\frac{1}{4},\frac{1}{3},\frac{5}{12},\frac{1}{2},\frac{7}{12},\frac{2}{3},\frac{3}{4},\frac{5}{6},\frac{11}{12},\frac{13}{12},\frac{7}{6};1\right)$$ $$=\frac{2}{3} \sqrt{2 \left(\sqrt{2}+1\right)}+\frac{\left(\frac{2 \Gamma \left(\frac{1}{4}\right)^2}{\Gamma \left(\frac{1}{12}\right)^2}+\frac{\sqrt[3]{2} \left(2 \sqrt{2}-\sqrt{3}+1\right) \Gamma \left(\frac{5}{24}\right)^2}{\Gamma \left(\frac{1}{24}\right)^2}\right) \Gamma \left(\frac{1}{3}\right)^2}{\pi }-\frac{2 \left(\sqrt{3}+1\right) \sqrt{\frac{2-\sqrt{2}}{\pi }} \Gamma \left(\frac{1}{4}\right) \Gamma \left(\frac{1}{3}\right)}{3^{5/8} \Gamma \left(\frac{1}{12}\right)}$$}\\
Meanwhile simplifying DIST($12,2,1,\{-\frac12\},\{\}$) gives (here $A=7 \sqrt{2}+16 \sqrt{3}+7 \sqrt{6}+8 \sqrt{3 \sqrt{3}+6}+28, B=-7 \sqrt{2}-16 \sqrt{3}+7 \sqrt{6}-8 \sqrt{6-3 \sqrt{3}}+28$):\\
{\footnotesize $$S=-\frac{\sqrt{A}}{6 \sqrt{2} \left(\sqrt{3}+2\right)^{3/4}}-\frac{\left(\sqrt{3}+2\right)^{3/4} }{\sqrt{6A}}+\frac{\left(2-\sqrt{3}\right)^{3/4} }{\sqrt{6B}}-\frac{\sqrt{B}}{6 \sqrt{2} \left(2-\sqrt{3}\right)^{3/4}}+\frac{\sqrt{2-\sqrt{2}} \sqrt[4]{2-\sqrt{3}}}{6 \sqrt{2}}$$ $$+\frac{\sqrt{\sqrt{2}+2} \sqrt[4]{2-\sqrt{3}}}{6 \sqrt{2}}+\frac{\sqrt{6-3 \sqrt{2}} \sqrt[4]{2-\sqrt{3}}}{6 \sqrt{2}}-\frac{\sqrt{\sqrt{2}+2} \sqrt[4]{2-\sqrt{3}}}{2 \sqrt{6}}-\frac{1}{6} \sqrt{\sqrt{2}+2} \sqrt[4]{\sqrt{3}+2}$$ $$-\frac{1}{6} \sqrt{6-3 \sqrt{2}} \sqrt[4]{\sqrt{3}+2}+\frac{\sqrt{2}}{\sqrt[4]{3}}+\frac{1}{6} \sqrt[4]{\sqrt{3}+2} \sqrt{2-\sqrt{2}}+\frac{\sqrt[4]{\sqrt{3}+2} \sqrt{2-\sqrt{2}}}{6 \sqrt{2}}+\frac{\sqrt{2-\sqrt{2}}}{3 \sqrt[4]{2}}+\frac{\sqrt[4]{\sqrt{3}+2} \sqrt{\sqrt{2}+2}}{6 \sqrt{2}}$$ $$+\frac{1}{3} \sqrt[4]{2} \sqrt{\sqrt{2}+2}+\frac{\sqrt{\sqrt{2}+2}}{3 \sqrt[4]{2}}+\frac{\sqrt[4]{\sqrt{3}+2} \sqrt{\sqrt{2}+2}}{2 \sqrt{6}}-\frac{2 \sqrt{2}}{3}-\frac{\sqrt{2}}{3^{3/4}}-\frac{\sqrt{6-3 \sqrt{2}} \sqrt[4]{\sqrt{3}+2}}{6 \sqrt{2}}+\frac{2}{\sqrt{3}}-\frac{\sqrt{\sqrt{2}+2} \sqrt[4]{\sqrt{3}+2}}{2 \sqrt{3}}$$}\\

\subsection{Polylog, K/E Gamma DIST}
Case $_3F_2$ (one may figure out the $E$-$O$ class of each result before DIST, modulo STIR/INIT/PFD):\\

$\text{DIST}\left(2,1,1,\left\{\frac{1}{2},1,1,\frac{5}{4},\frac{4}{3}\right\},\left\{\frac{1}{4},\frac{1}{3},\frac{5}{2},\frac{5}{2}\right\}\right)$\\
$$\, _7F_6\left(\frac{3}{4},1,1,\frac{5}{4},\frac{3}{2},\frac{13}{8},\frac{5}{3};\frac{5}{8},\frac{2}{3},\frac{7}{4},\frac{7}{4},\frac{9}{4},\frac{9}{4};1\right)=\frac{3195 C}{64}+\frac{3105 \text{Li}_2\left(1-\sqrt{2}\right)}{128}$$ $$-\frac{3105 \text{Li}_2\left(\sqrt{2}-1\right)}{128}+\frac{3105 \pi ^2}{512}-\frac{5265}{64}-\frac{1}{128} 3105 \log ^2\left(\sqrt{2}+1\right)+\frac{2115 \log \left(\sqrt{2}+1\right)}{64 \sqrt{2}}$$\\

$\text{DIST}\left(2,0,1,\left\{\frac{1}{2},\frac{3}{2},2,\frac{5}{2}\right\},\left\{\frac{7}{2},\frac{7}{2},4\right\}\right)$\\
$$\, _8F_7\left(\frac{1}{4},\frac{3}{4},\frac{3}{4},1,\frac{5}{4},\frac{5}{4},\frac{3}{2},\frac{7}{4};\frac{1}{2},\frac{7}{4},\frac{7}{4},2,\frac{9}{4},\frac{9}{4},\frac{5}{2};1\right)=-\frac{675 \text{Li}_2\left(1-\sqrt{2}\right)}{4}+\frac{675}{4} \text{Li}_2\left(-\frac{1}{\sqrt{2}}\right)$$ $$+\frac{225 \pi ^2}{16}+\frac{2995}{8 \sqrt{2}}+\frac{375 \pi }{32}-600+\frac{675 \log ^2(2)}{32}+\frac{675}{8} \log \left(\sqrt{2}+1\right) \log (2)+\frac{675}{8} \pi  \log (2)-\frac{825}{16} \log \left(\sqrt{2}+1\right)$$\\

$\text{DIST}\left(2,0,1,\left\{\frac{1}{2},\frac{1}{2},1,\frac{4}{3},\frac{4}{3}\right\},\left\{\frac{1}{3},\frac{1}{3},3,3\right\}\right)$\\
$$\, _7F_6\left(\frac{1}{4},\frac{1}{4},1,\frac{3}{4},\frac{3}{4},\frac{7}{6},\frac{7}{6};\frac{1}{6},\frac{1}{6},\frac{3}{2},\frac{3}{2},2,2;1\right)=-\frac{1600}{9}+\frac{10336}{27 \pi }-\frac{80 \sqrt{2} \Gamma \left(\frac{1}{4}\right)^2}{27 \pi ^{3/2}}+\frac{352 \sqrt{2 \pi }}{\Gamma \left(\frac{1}{4}\right)^2}$$\\

$\text{DIST}\left(2,1,1,\left\{\frac{1}{2},1,\frac{7}{6},\frac{4}{3},\frac{3}{2}\right\},\left\{\frac{1}{6},\frac{1}{3},2,4\right\}\right)$\\
$$\, _7F_6\left(1,\frac{3}{4},\frac{5}{4},\frac{5}{4},\frac{19}{12},\frac{5}{3},\frac{7}{4};\frac{7}{12},\frac{2}{3},\frac{3}{2},\frac{5}{2},2,3;1\right)=\frac{320}{7}-\frac{19328}{315 \pi }-\frac{2224 \sqrt{2} \Gamma \left(\frac{1}{4}\right)^2}{63 \pi ^{3/2}}+\frac{52352 \sqrt{2 \pi }}{105 \Gamma \left(\frac{1}{4}\right)^2}$$\\

\noindent Case $_4F_3$:\\

$\text{DIST}\left(2,0,1,\left\{\frac{1}{2},\frac{1}{2},\frac{1}{2},1,\frac{4}{3},\frac{4}{3}\right\},\left\{\frac{1}{3},\frac{1}{3},2,2,2\right\}\right)$\\
$$\, _8F_7\left(\frac{1}{4},\frac{1}{4},\frac{1}{4},\frac{3}{4},\frac{3}{4},\frac{3}{4},\frac{7}{6},\frac{7}{6};\frac{1}{6},\frac{1}{6},1,1,\frac{3}{2},\frac{3}{2},\frac{3}{2};1\right)$$ $$=\frac{15}{\pi }-\frac{15 \Gamma \left(\frac{1}{4}\right)^4}{8 \pi ^3}+\frac{35 \Gamma \left(\frac{1}{8}\right)^2 \Gamma \left(\frac{3}{8}\right)^2}{32 \sqrt{2} \pi ^3}+\frac{264 \pi  \sqrt{2}}{\Gamma \left(\frac{1}{8}\right)^2 \Gamma \left(\frac{3}{8}\right)^2}-\frac{264 \pi }{\Gamma \left(\frac{1}{4}\right)^4}$$\\

$\text{DIST}\left(2,1,1,\left\{\frac{1}{2},\frac{1}{2},\frac{1}{2},1,\frac{5}{4},\frac{7}{4}\right\},\left\{\frac{1}{4},\frac{3}{4},2,2,3\right\}\right)$\\
$$\, _9F_8\left(\frac{3}{4},\frac{3}{4},\frac{3}{4},1,\frac{5}{4},\frac{5}{4},\frac{5}{4},\frac{13}{8},\frac{15}{8};\frac{5}{8},\frac{7}{8},\frac{3}{2},\frac{3}{2},2,2,2,\frac{5}{2};1\right)$$ $$=\frac{4608}{35}+\frac{512}{105 \pi }-\frac{4288 \Gamma \left(\frac{1}{4}\right)^4}{315 \pi ^3}-\frac{24 \sqrt{2} \Gamma \left(\frac{1}{8}\right)^2 \Gamma \left(\frac{3}{8}\right)^2}{5 \pi ^3}+\frac{4096 \pi }{5 \Gamma \left(\frac{1}{4}\right)^4}-\frac{4096 \sqrt{2} \pi }{45 \Gamma \left(\frac{1}{8}\right)^2 \Gamma \left(\frac{3}{8}\right)^2}$$\\

$\text{DIST}\left(2,0,1,\left\{\frac{1}{2},\frac{1}{2},\frac{1}{2},\frac{1}{2}\right\},\left\{\frac{3}{2},\frac{3}{2},\frac{3}{2}\right\}\right)$\\
$$\, _5F_4\left(\frac{1}{4},\frac{1}{4},\frac{1}{4},\frac{1}{4},\frac{3}{4};\frac{1}{2},\frac{5}{4},\frac{5}{4},\frac{5}{4};1\right)=-\frac{1}{8} \text{Li}_3\left(\frac{1}{2}-\frac{1}{\sqrt{2}}\right)+\frac{1}{8} \text{Li}_3\left(2 \sqrt{2}-2\right)+\frac{\pi ^3}{96}+\frac{\log ^3(2)}{48}-\frac{1}{16} \log ^3\left(\sqrt{2}+1\right)$$ $$+\frac{1}{16} \log \left(\sqrt{2}+1\right) \log ^2(2)+\frac{1}{8} \pi  \log ^2(2)+\frac{1}{16} \log ^2\left(\sqrt{2}+1\right) \log (2)+\frac{1}{48} \pi ^2 \log (2)+\frac{1}{48} \pi ^2 \log \left(\sqrt{2}+1\right)$$\\

$\text{DIST}\left(2,0,1,\left\{1,1,1,\frac{3}{2},\frac{5}{2}\right\},\left\{\frac{1}{2},3,3,4\right\}\right)$\\
$$\, _8F_7\left(\frac{1}{2},\frac{1}{2},1,1,1,\frac{5}{4},\frac{5}{4},\frac{7}{4};\frac{1}{4},\frac{3}{2},\frac{3}{2},2,2,2,\frac{5}{2};1\right)=-88 \text{Li}_2\left(\frac{1}{2}-\frac{1}{\sqrt{2}}\right)-84 \sqrt{2}-\frac{26 \pi ^2}{3}+128$$ $$+148 \log ^2(2)+44 \log ^2\left(\sqrt{2}+1\right)-88 \log \left(\sqrt{2}+1\right) \log (2)-64 \log (2)+60 \log \left(\sqrt{2}+1\right)$$\\

$\text{DIST}\left(2,0,1,\left\{1,1,1,1,\frac{5}{4},\frac{5}{2}\right\},\left\{\frac{1}{4},\frac{3}{2},\frac{3}{2},3,3\right\}\right)$\\
$$\, _9F_8\left(\frac{1}{2},\frac{1}{2},\frac{1}{2},1,1,1,1,\frac{9}{8},\frac{7}{4};\frac{1}{8},\frac{3}{4},\frac{3}{4},\frac{5}{4},\frac{3}{2},\frac{3}{2},2,2;1\right)=\frac{13}{3} \text{Li}_3\left(3-2 \sqrt{2}\right)+\frac{26}{3} \text{Li}_2\left(3-2 \sqrt{2}\right) \log \left(\sqrt{2}+1\right)$$ $$-\frac{59 \zeta (3)}{12}-\frac{\pi ^2}{4}+\frac{52}{9} \log ^3\left(\sqrt{2}+1\right)+\frac{19}{3} \log ^2\left(\sqrt{2}+1\right)-\frac{26}{3} \log (2) \log ^2\left(\sqrt{2}+1\right)+2 \sqrt{2} \log \left(\sqrt{2}+1\right)+\frac{1}{6} \pi ^2 \log (2)$$\\

\clearpage

\noindent Case $_5F_4$:\\

$\text{DIST}\left(2,1,1,\left\{1,1,1,1,\frac{3}{2}\right\},\{2,2,2,2\}\right)$\\
$$\, _6F_5\left(1,1,1,1,\frac{5}{4},\frac{7}{4};\frac{3}{2},2,2,2,2;1\right)=-\frac{64}{3} \text{Li}_3\left(2 \sqrt{2}-2\right)+\frac{64}{3} \text{Li}_3\left(\frac{1}{2}-\frac{1}{\sqrt{2}}\right)-\frac{64}{3} \text{Li}_2\left(2 \sqrt{2}-2\right) \log \left(\sqrt{2}+1\right)$$ $$+\frac{64}{3} \text{Li}_2\left(\frac{1}{2}-\frac{1}{\sqrt{2}}\right) \log \left(\sqrt{2}+1\right)+\frac{128 \zeta (3)}{3}+\frac{224 \log ^3(2)}{9}-\frac{64}{3} \log ^3\left(\sqrt{2}+1\right)+\frac{32}{3} \log ^2\left(\sqrt{2}+1\right) \log (2)$$\\

$\text{DIST}\left(2,0,1,\{1,1,1,1,1\},\left\{\frac{3}{2},2,2,2\right\}\right)$\\
{\footnotesize $$\, _6F_5\left(\frac{1}{2},\frac{1}{2},\frac{1}{2},\frac{1}{2},1,1;\frac{3}{4},\frac{5}{4},\frac{3}{2},\frac{3}{2},\frac{3}{2};1\right)=2 \text{Li}_4\left(\frac{1}{2}\right)-\text{Li}_4\left(\frac{1}{2}-\frac{1}{\sqrt{2}}\right)+\text{Li}_4\left(2 \sqrt{2}-2\right)-\frac{1}{2} \text{Li}_2\left(\frac{1}{2}-\frac{1}{\sqrt{2}}\right) \log ^2(2)$$ $$+\frac{1}{2} \text{Li}_2\left(2 \sqrt{2}-2\right) \log ^2(2)-\frac{1}{2} \text{Li}_2\left(\frac{1}{2}-\frac{1}{\sqrt{2}}\right) \log ^2\left(\sqrt{2}+1\right)+\frac{5}{2} \text{Li}_2\left(2 \sqrt{2}-2\right) \log ^2\left(\sqrt{2}+1\right)$$ $$+2 \text{Li}_2\left(3-2 \sqrt{2}\right) \log ^2\left(\sqrt{2}+1\right)-\text{Li}_2\left(\frac{1}{2}-\frac{1}{\sqrt{2}}\right) \log \left(\sqrt{2}+1\right) \log (2)-\text{Li}_2\left(2 \sqrt{2}-2\right) \log \left(\sqrt{2}+1\right) \log (2)$$ $$-2 \text{Li}_2\left(3-2 \sqrt{2}\right) \log \left(\sqrt{2}+1\right) \log (2)-\text{Li}_3\left(\frac{1}{2}-\frac{1}{\sqrt{2}}\right) \log (2)-\text{Li}_3\left(2 \sqrt{2}-2\right) \log (2)-\text{Li}_3\left(3-2 \sqrt{2}\right) \log (2)$$ $$-\text{Li}_3\left(\frac{1}{2}-\frac{1}{\sqrt{2}}\right) \log \left(\sqrt{2}+1\right)+\text{Li}_3\left(2 \sqrt{2}-2\right) \log \left(\sqrt{2}+1\right)+\zeta (3) \log (2)-\frac{47 \pi ^4}{1440}-\frac{1}{24} \log ^4(2)+\frac{35}{8} \log ^4\left(\sqrt{2}+1\right)$$ $$-\frac{1}{2} \log \left(\sqrt{2}+1\right) \log ^3(2)-\frac{35}{6} \log ^3\left(\sqrt{2}+1\right) \log (2)+\frac{1}{6} \pi ^2 \log ^2(2)+\frac{9}{4} \log ^2\left(\sqrt{2}+1\right) \log ^2(2)-\frac{1}{3} \pi ^2 \log ^2\left(\sqrt{2}+1\right)$$}\\

$\text{DIST}\left(2,0,\frac{1}{64},\left\{1,1,1,1,\frac{3}{2}\right\},\{2,2,2,2\}\right)$\\
$$\, _5F_4\left(\frac{1}{2},\frac{1}{2},\frac{1}{2},1,1;\frac{3}{4},\frac{5}{4},\frac{3}{2},\frac{3}{2};\frac{1}{64}\right)=-4 \pi  \Im\left(\text{Li}_2\left(\frac{3}{4}-\frac{i \sqrt{7}}{4}\right)\right)+8 \tan ^{-1}\left(\sqrt{7}\right) \Im\left(\text{Li}_2\left(\frac{3}{4}-\frac{i \sqrt{7}}{4}\right)\right)$$ $$+4 \Re\left(\text{Li}_3\left(\frac{3}{4}-\frac{i \sqrt{7}}{4}\right)\right)-\frac{7 \zeta (3)}{2}-\frac{1}{3} \log ^3(2)-\pi ^2 \log (2)-4 \log (2) \tan ^{-1}\left(\sqrt{7}\right)^2+4 \pi  \log (2) \tan ^{-1}\left(\sqrt{7}\right)$$\\

$\text{DIST}\left(2,0,\frac{1}{16},\{1,1,1,1,1\},\left\{\frac{3}{2},2,2,2\right\}\right)$\\
$$\, _6F_5\left(\frac{1}{2},\frac{1}{2},\frac{1}{2},\frac{1}{2},1,1;\frac{3}{4},\frac{5}{4},\frac{3}{2},\frac{3}{2},\frac{3}{2};\frac{1}{16}\right)=8 \text{Li}_4\left(\frac{\sqrt{5}}{2}-\frac{1}{2}\right)-\frac{1}{2} \text{Li}_4\left(\frac{3}{2}-\frac{\sqrt{5}}{2}\right)+4 \text{Li}_3\left(\frac{\sqrt{5}}{2}-\frac{1}{2}\right) \log \left(\frac{\sqrt{5}}{2}+\frac{1}{2}\right)$$ $$-4 \text{Li}_3\left(\frac{1}{2}-\frac{\sqrt{5}}{2}\right) \log \left(\frac{\sqrt{5}}{2}+\frac{1}{2}\right)-\frac{47 \pi ^4}{648}-\frac{3}{2}  \log ^4\left(\frac{\sqrt{5}}{2}+\frac{1}{2}\right)+\frac{1}{3} \pi ^2 \log ^2\left(\frac{\sqrt{5}}{2}+\frac{1}{2}\right)$$\\

\noindent Bonus: 2 more examples generated by DIST (and STIR $\cdots$) on certain functions.\\

$$\, _6F_5\left(\frac{1}{4},\frac{1}{2},\frac{1}{2},\frac{1}{2},1,1;\frac{3}{4},\frac{5}{4},\frac{5}{4},\frac{3}{2},\frac{3}{2};1\right)=4 C+\frac{1}{2} \text{Li}_3\left(3-2 \sqrt{2}\right)+\text{Li}_2\left(3-2 \sqrt{2}\right) \log \left(\sqrt{2}+1\right)$$ $$+\frac{3 \zeta (3)}{8}+\frac{2}{3} \log ^3\left(\sqrt{2}+1\right)-\log (2) \log ^2\left(\sqrt{2}+1\right)-2 \log ^2\left(\sqrt{2}+1\right)-\frac{1}{4} \pi ^2 \log (2)$$\\

$$\, _6F_5\left(\frac{1}{4},\frac{1}{4},\frac{3}{4},\frac{3}{4},\frac{3}{4},\frac{5}{4};\frac{1}{2},\frac{7}{4},\frac{7}{4},\frac{9}{4},\frac{9}{4};1\right)=-\frac{1125 \text{Li}_2\left(1-\sqrt{2}\right)}{256}+\frac{1125}{256} \text{Li}_2\left(-\frac{1}{\sqrt{2}}\right)+\frac{375 \pi ^2}{1024}+\frac{225}{64 \sqrt{2}}$$ $$+\frac{225 \pi }{64}+\frac{1125 \log ^2(2)}{2048}+\frac{1125}{512} \log \left(\sqrt{2}+1\right) \log (2)-\frac{2475}{512} \pi  \log (2)-\frac{225}{32} \log \left(\sqrt{2}+1\right)$$\\
One may manually figure out the original function for these low-order cases, but\\

\noindent Note 4: An algorithmic procedure detecting DIST reducibility of high order $_pF_q$ is inaccessible in general since there exists hidden elimination between parameters.\\

\noindent Note 5: Polylog arguments in above closed-forms are not completely regularized, thus may be furtherly simplified using certain functional equations. Those $_pF_q$ expressible via polylogs with argument involving $\sqrt{2}$ have close relationship with level 8 MZVs (multiple Zeta values). For a general statement of corresponding level 4 one may see \cite{ZMH}, but we are not aware of a generalization to level 8 here since method of DIST is rather limited compared to that article.\\

\noindent Note 6: The one given by $\text{DIST}\left(2,0,\frac{1}{16},\{1,1,1,1,1\},\left\{\frac{3}{2},2,2,2\right\}\right)$ is related to level 5 MZVs. One may replace $z=\frac{1}{16}$ by $z=\frac{1}{4}$ to obtain a similar result featuring level 4 and 12. What about $z=\frac1{64}$?\\

\subsection{Miscellaneous DIST}
\noindent 1. Using trigonometric-hyperbolic representation of Bessel $_0F_1$ at half-integer parameters, one may evaluate some Mittag-Leffler type series. For example DIST($12,0,1,\{\},\{\frac12\}$)  gives:\\
$$\, _0F_{23}\left(/;\frac{1}{24},\frac{1}{12},\frac{1}{8},\frac{1}{6},\frac{5}{24},\frac{1}{4},\frac{7}{24},\frac{1}{3},\frac{3}{8},\frac{5}{12},\frac{11}{24},\frac{1}{2},\frac{13}{24},\frac{7}{12},\frac{5}{8},\frac{2}{3},\frac{17}{24},\frac{3}{4},\frac{19}{24},\frac{5}{6},\frac{7}{8},\frac{11}{12},\frac{23}{24};1\right)$$ $$=\frac{\cos (24)}{12}+\frac{\cosh (24)}{12}+\frac{1}{6} \cos \left(12 \sqrt{3}\right) \cosh (12)+\frac{1}{6} \cos \left(12 \sqrt{2}\right) \cosh \left(12 \sqrt{2}\right)+\frac{1}{6} \cos (12) \cosh \left(12 \sqrt{3}\right)$$ $$+\frac{1}{6} \cos \left(6 \sqrt{2}-6 \sqrt{6}\right) \cosh \left(6 \sqrt{2}+6 \sqrt{6}\right)+\frac{1}{6} \cos \left(6 \sqrt{2}+6 \sqrt{6}\right) \cosh \left(6 \sqrt{2}-6 \sqrt{6}\right)$$\\
This can be generalized to arbitrary order like in Note 3. One may also try $_0F_2$ and etc which we omit.\\

\noindent 2. All $\sum_{n=0}^{\infty} {\binom{2n}n}^{\pm1}n^m z^n$ are expressible in elementary functions by differentiation. One may write them in terms of $_pF_q$ and take the DIST procedure, for instance:\\
$$\, _6F_5\left(\frac{5}{4},\frac{7}{4},2,2,2,2;1,1,1,1,\frac{3}{2};\frac{1}{4}\right)=\frac{31}{2592 \sqrt{6}}+\frac{4921}{96 \sqrt{2}}$$

$$\, _6F_5\left(\frac{3}{2},2,2,2,2,2;1,1,1,\frac{5}{4},\frac{7}{4};\frac{1}{64}\right)=\frac{617416}{583443}+\frac{19 \log (2)}{729}+\frac{5298 \cot ^{-1}\left(\sqrt{7}\right)}{2401 \sqrt{7}}$$\\
This can be generalized to $\sum_{n=0}^{\infty} {\binom{3n}n}^{\pm1}n^m z^n$ using $\sum _{n=0}^{\infty } \binom{3 n}{n} x^n=\frac{2 \cos \left(\frac{1}{6} \cos ^{-1}\left(1-\frac{27 x}{2}\right)\right)}{\sqrt{4-27 x}}$ (related to Lagrange inversion) and the complicated log-algebraic closed-form of the inverse case given in \cite{inv}. Now taking DIST again and plug in special values, for instance:\\
$$\, _4F_3\left(\frac{1}{6},\frac{1}{3},\frac{2}{3},\frac{5}{6};\frac{1}{4},\frac{1}{2},\frac{3}{4};\frac{729}{1024}\right)=\frac{2 \sqrt[3]{\frac{2}{3 \sqrt{3}+\sqrt{59}}}}{\sqrt{59}}+\frac{\sqrt[3]{3 \sqrt{6}+\sqrt{118}}}{\sqrt{118}}+2 \sqrt{\frac{2}{5} \left(\frac{3 \sqrt{5}}{16}+\frac{7}{16}\right)}$$\\

\noindent 3. This part is originally due to \cite{sos}. Consider evaluating $\int_0^{\frac{\pi}2}\tan^{-1}(\sin^3(x))dx$ in 2 ways: One is to power-expand $\tan^{-1}$ and apply Wallis formula (yielding hypergeometric form), another is to factorize the integrand into $\tan^{-1}(w^k z)\ (w=e^{\frac{2 \pi i}3})$ and apply Fourier series of $\tan^{-1}(r\sin(x))$ (yielding polylog form). Equating these 2 forms gives\\
$$\, _5F_4\left(\frac{1}{2},\frac{2}{3},1,1,\frac{4}{3};\frac{5}{6},\frac{7}{6},\frac{3}{2},\frac{3}{2};-1\right)=\frac{3 \pi ^2}{32}+\frac{3}{4} \log ^2\left(\sqrt{2}+1\right)-\frac{3}{8} \log ^2\left(2-\sqrt{3}\right)$$\\
This method, essentially a order 3 DIST of $_3F_2$ (find it), can be generalized to $\tan^{-1}(\sin^{2k+1}(x))$ for arbitrary $k\in \mathbb N$. The specific case, however, is of own interest since the closed-form is related to level 8 and 12 MZVs, suggesting deeper relations between them and $_pF_q$ having parameters of denominator 6. What about results generated by corresponding $\log(1+\sin^k(x))$?\\

\noindent 4.  Some $_2F_1$ and $_3F_2$ have elliptic closed-forms due to quadratic transforms, Clausen formula and contiguous relations. For instance $\text{DIST}\left(2,0,1,\left\{\frac{1}{4},\frac{1}{2},\frac{3}{4}\right\},\{1,1\}\right)$ gives:\\
$$\, _6F_5\left(\frac{1}{8},\frac{1}{4},\frac{3}{8},\frac{5}{8},\frac{3}{4},\frac{7}{8};\frac{1}{2},\frac{1}{2},\frac{1}{2},1,1;1\right)=\frac{32 \pi  \sqrt{2 \sqrt{2}-2} K\left(\frac{1}{2}-\sqrt{\frac{1}{\sqrt{2}}-\frac{1}{2}}\right)^2+\Gamma \left(\frac{1}{8}\right)^2 \Gamma \left(\frac{3}{8}\right)^2}{16 \pi ^3}$$\\
One may choose appropriate arguments to make use of more elliptic singular values. For example $\text{DIST}\left(2,0,\frac{1}{4096},\left\{\frac{3}{2},\frac{3}{2},\frac{3}{2}\right\},\{2,2\}\right)$ and plugging in  $K(k_7)$ gives:\\
$$\, _6F_5\left(\frac{3}{4},\frac{3}{4},\frac{3}{4},\frac{5}{4},\frac{5}{4},\frac{5}{4};\frac{1}{2},1,1,\frac{3}{2},\frac{3}{2};\frac{1}{4096}\right)=\frac{512 K\left(\frac{1}{2}-\frac{\sqrt{65}}{16}\right)^2}{\pi ^2}+\frac{4096 K\left(\frac{1}{2}-\frac{\sqrt{65}}{16}\right)^2}{\sqrt{65} \pi ^2}$$ $$-\frac{8192 K\left(\frac{1}{2}-\frac{\sqrt{65}}{16}\right) E\left(\frac{1}{2}-\frac{\sqrt{65}}{16}\right)}{\sqrt{65} \pi ^2}+\frac{2048}{21 \pi }-\frac{160 \Gamma \left(\frac{1}{7}\right)^2 \Gamma \left(\frac{2}{7}\right)^2 \Gamma \left(\frac{4}{7}\right)^2}{21 \sqrt{7} \pi ^4}$$\\
The case replacing $\frac1{4096}$ by $\frac 14$ is also interesting due to $K(k_3)$. For $_4F_3$ we have an isolated $K$ formula (originally due to \cite{Bry}, which has a typo):\\
$$\, _4F_3\left(\frac{1}{4},\frac{1}{4},\frac{3}{4},\frac{3}{4};\frac{1}{2},1,1;z\right)=\frac{4}{\pi ^2} K\left(\frac{1}{2}-\frac{1}{2 \sqrt{-2 z-2 \sqrt{z^2-z}+1}}\right) K\left(\frac{1}{2}-\frac{1}{2} \sqrt{-2 z-2 \sqrt{z^2-z}+1}\right)$$\\
Which can be established by $z\to-\left(\frac t2 - \frac 1{2t}\right)^2$, an argument similar to proof of Lemma 3 (use `DifferentialRootReduce' to check) and analytic continuation. Thus $\text{DIST}\left(2,0,1,\left\{\frac{1}{4},\frac{1}{4},\frac{3}{4},\frac{3}{4}\right\},\left\{\frac{1}{2},1,1\right\}\right)$ yields:\\
$$\, _8F_7\left(\frac{1}{8},\frac{1}{8},\frac{3}{8},\frac{3}{8},\frac{5}{8},\frac{5}{8},\frac{7}{8},\frac{7}{8};\frac{1}{4},\frac{1}{2},\frac{1}{2},\frac{1}{2},\frac{3}{4},1,1;1\right)=\frac{2 \left(K\left(1-\frac{1}{\sqrt{2}}\right) K\left(-\frac{1}{\sqrt{2}}\right)+K\left(\frac{1}{2}+\frac{i}{2}\right) K\left(\frac{1}{2}-\frac{i}{2}\right)\right)}{\pi ^2}$$\\

\noindent 5. One may evaluate similar hypergeometric series in polylog closed-forms using various techniques other than DIST relations. Here we illustrate some examples arising from section `Denominator 4'. Firstly integrating $\, _2F_1\left(\frac{5}{4},\frac{7}{4};\frac{3}{2};z\right)=\frac{\frac{1}{\left(1-\sqrt{z}\right)^{3/2}}-\frac{1}{\left(\sqrt{z}+1\right)^{3/2}}}{3 \sqrt{z}}$ 2 times one may compute $\, _4F_3\left(1,1,\frac{5}{4},\frac{7}{4};\frac{3}{2},2,2;z\right)$ in closed forms. Multiply it with Beta kernel $\frac{1}{\sqrt{1-z}}$ and integrate on $(0,1)$ directly yields:\\
$$\, _5F_4\left(1,1,1,\frac{5}{4},\frac{7}{4};\frac{3}{2},\frac{3}{2},2,2;1\right)=-\frac{8}{3} \text{Li}_2\left(3-2 \sqrt{2}\right)+\frac{2 \pi ^2}{9}-\frac{8}{3}  \log ^2(2)-\frac{8}{3} \log ^2\left(\sqrt{2}+1\right)+\frac{16}{3} \log \left(\sqrt{2}+1\right) \log (2)$$\\
By similar means, integrating $\, _2F_1\left(\frac{1}{4},\frac{3}{4};\frac{1}{2};z\right)=\frac{\sqrt{\sqrt{1-z}+1}}{\sqrt{2-2 z}}$ 3 times indefinitely by brute-force and letting $z\to1^-$ give closed-forms of $\, _4F_3\left(\frac{1}{4},\frac{3}{4},1,1;\frac{1}{2},2,2;1\right), \, _5F_4\left(\frac{1}{4},\frac{3}{4},1,1,1;\frac{1}{2},2,2,2;1\right)$, which make the following constructible via PFD and elementary terms:\\
$$\, _7F_6\left(\frac{1}{4},\frac{1}{4},\frac{3}{4},\frac{3}{4},1,1,1;\frac{1}{2},\frac{5}{4},\frac{7}{4},2,2,2;1\right)=-\frac{8}{3} \text{Li}_2\left(\frac{2-\sqrt{2}}{4} \right)+\frac{8}{3} \text{Li}_2\left(\frac{2+\sqrt{2}}{4} \right)+\frac{32 \text{Li}_2\left(\frac{1}{\sqrt{2}}\right)}{3}$$ $$-2 \pi ^2-\frac{664 \sqrt{2}}{9}-\frac{100 \pi }{9}+\frac{1280}{9}+\frac{68 \log ^2(2)}{3}-8 \log \left(\sqrt{2}+1\right) \log (2)-\frac{832 \log (2)}{9}+72 \log \left(\sqrt{2}+1\right)$$\\
If we preserve $\, _5F_4\left(\frac{1}{4},\frac{3}{4},1,1,1;\frac{1}{2},2,2,2;z\right)$ obtained in above process and compute the corresponding $\text{DIST}\left(2,1,1,\left\{\frac{1}{4},\frac{3}{4},1,1,1\right\},\left\{\frac{1}{2},2,2,2\right\}\right)$, one arrive at the following after simplifications (as one may observe it seems to be related to level 16 MZVs):\\

{\footnotesize $$\, _7F_6\left(\frac{5}{8},\frac{7}{8},1,1,1,\frac{9}{8},\frac{11}{8};\frac{3}{4},\frac{5}{4},\frac{3}{2},2,2,2;1\right)=-\frac{256}{9} \text{Li}_2\left(\frac{1}{4} \left(2-\sqrt{2}\right)\right)+\frac{256}{9} \text{Li}_2\left(\frac{1}{4} \left(\sqrt{2}+2\right)\right)$$ $$-\frac{256}{9} \text{Li}_2\left(\frac{1}{4} \left(2-\sqrt{2 \sqrt{2}+2}\right)\right)+\frac{1024 \text{Li}_2\left(\frac{1}{\sqrt{2}}\right)}{9}+\frac{512}{9} \text{Li}_2\left(-\sqrt{2 \sqrt{2}-2}\right)-\frac{512}{9} \text{Li}_2\left(\sqrt{2 \sqrt{2}-2}\right)$$ $$-\frac{256}{9} \text{Li}_2\left(\frac{4}{\sqrt{2 \sqrt{2}+2}+2}\right)-\frac{13312}{81} \sqrt{2 \left(\sqrt{2}+1\right)}-\frac{64 \pi ^2}{27}-\frac{13312 \sqrt{2}}{81}-\frac{2048 \sqrt{\sqrt{2}+1}}{81}+\frac{57344}{81}+\frac{1312 \log ^2(2)}{3}$$ $$-\frac{128}{9} \log ^2\left(\sqrt{2}+1\right)-\frac{128}{9} \log ^2\left(\sqrt{2}+\sqrt{\sqrt{2}+1}\right)-128 \log \left(\sqrt{2}+1\right) \log (2)-\frac{128}{3} \log \left(\sqrt{2}+\sqrt{\sqrt{2}+1}\right) \log (2)$$ $$-\frac{20480 \log (2)}{27}+\frac{2560}{9} \log \left(\sqrt{2}+1\right)+\frac{5120}{27} \log \left(\sqrt{2}+\sqrt{\sqrt{2}+1}\right)-\frac{256}{9} \log \left(\sqrt{2}+1\right) \log \left(\sqrt{2}+\sqrt{\sqrt{2}+1}\right)$$}\\
A final example: By Pfaff transform one have $\, _2F_1\left(\frac{3}{4},\frac{3}{4};\frac{7}{4};z\right)=-\frac{3 \left(\tan ^{-1}\left(\sqrt[4]{\frac{z}{z-1}}\right)-\tanh ^{-1}\left(\sqrt[4]{\frac{z}{z-1}}\right)\right)}{2 (-z)^{3/4}}$, thus $z^{-\frac34}\int \frac{\, _2F_1\left(\frac{3}{4},\frac{3}{4};\frac{7}{4};z\right)}{\sqrt[4]{z}} \, dz$ can be calculated using $z\to \frac{t^4}{1+t^4}$ and brute force, which leads to closed-form of $\, _3F_2\left(\frac{3}{4},\frac{3}{4},\frac{3}{4};\frac{7}{4},\frac{7}{4};z\right)$ up to constant. Setting $z=\frac12$ and simplifying lead to a result analogous to the problem of \cite{Vla}, whose author also suggests the closed-form of $\, _4F_3\left(\frac{1}{2},\frac{1}{2},\frac{1}{2},\frac{1}{2};\frac{3}{2},\frac{3}{2},\frac{3}{2};z\right)$ in section 3.1:\\
$$2^{3/4} \, _3F_2\left(\frac{3}{4},\frac{3}{4},\frac{3}{4};\frac{7}{4},\frac{7}{4};\frac{1}{2}\right)=\left(\frac{9}{2}+\frac{9 i}{2}\right) \text{Li}_2\left(-i \left(\sqrt{2}-1\right)\right)+\left(\frac{9}{2}-\frac{9 i}{2}\right) \text{Li}_2\left(i \left(\sqrt{2}-1\right)\right)$$ $$-\frac{9}{8} \text{Li}_2\left(3-2 \sqrt{2}\right)-\frac{3 \pi ^2}{4}+\frac{9}{8} \log ^2\left(\sqrt{2}+1\right)-\frac{9}{2} \log (2) \log \left(\sqrt{2}+1\right)+\frac{9}{8} \pi  \log \left(\sqrt{2}+1\right)+\frac{9}{4} \pi  \log (2)$$\\
This is related to level 8 again. All these calculations and simplifications are done with help of self-written Mathematica commands.\\\\

\noindent \textbf{\large Acknowledgements.} Special thanks to author of \cite{sos} and \cite{Vla} for their heuristic results.\\

\end{document}